\documentclass[12pt, a4paper]{amsart}

\usepackage{amsmath,amsthm,amsfonts,mathrsfs}
\usepackage{latexsym}
\usepackage[abbrev]{amsrefs}
\usepackage{bm}
\usepackage{url}
\usepackage[psamsfonts]{amssymb}
\usepackage{ascmac}

\theoremstyle{definition}
\theoremstyle{plain}
 \newtheorem{thm}{Theorem}[section]
 \newtheorem{lem}[thm]{Lemma}
 \newtheorem{cor}[thm]{Corollary}
 \newtheorem{prop}[thm]{Proposition}
 \newtheorem{claim}[thm]{Claim}
\theoremstyle{definition}
 \newtheorem{defi}[thm]{Definition}
\theoremstyle{remark}
 \newtheorem{rem}[thm]{\it Remark}

\newtheorem*{Acknowledgements}{\it Acknowledgements}

\newenvironment{pf}{\par\begin{trivlist}%
\item[]{\bf Proof.}\ }{\hfill $\square$ \end{trivlist}\par}
\newenvironment{apf}[1]{\par\begin{trivlist}%
\item[]{\bf Proof of #1.}\ }{\hfill $\square$ \end{trivlist}\par}
\newenvironment{cpf}{\par\begin{trivlist}%
\item[]{\bf Proof.}\ }{\hfill $\Diamond$ \end{trivlist}\par}

\newcommand{\Ric}{\mathrm{Ric}}
\newcommand{\diam}{\mathrm{diam}}
\newcommand{\vol}{\mathrm{vol}}
\newcommand{\vb}{v_{cK,b}}
\newcommand{\va}{v_{cK,a}}
\newcommand{\hess}{\mathrm{Hess}}
\newcommand{\sn}{{\bf sn}}
\newcommand{\vep}{\varepsilon}
\newcommand{\lmd}{\lambda}
\newcommand{\R}{\mathbb{R}}
\newcommand{\N}{\mathbb{N}}

\newcommand{\A}{{\mathcal{A}}}

\theoremstyle{defi}

\def\lf{\left}
\def\r{\right}

\BibSpec{arXiv}{%
  +{}{\PrintAuthors}{author}
  +{,}{ \textit}{title}
  +{}{ \parenthesize}{date}
  +{,}{ arXiv }{eprint}
}

\numberwithin{equation}{section}

\begin{document}
\title[]{Quantitative estimate of diameter \\for weighted manifolds \\under integral curvature bounds and $\vep$-range}
\author[]{Taku Ito}
\address{\hspace{-10pt}\scriptsize$\ast$ Department of Mathematical Sciences, Tokyo Metropolitan University, Tokyo {192-0397}, Japan.}
\email{size14sphere@gmail.com}
\date{\today}
\vspace{-10pt}
\begin{abstract}
In this article, we extend the compactness theorems proved by Sprouse~\cite{SC} and Hwang--Lee~\cite{HL} to a weighted manifold under the assumption that the weighted Ricci curvature is bounded below in terms of its weight function. With the help of the $\vep$-range, we treat the case that the effective dimension is at most $1$ in addition to the case that the effective dimension is at least the dimension of the manifold. To show these theorems, we extend the segment inequality of Cheeger--Colding~\cite{CC} to a weighted manifold. 
\end{abstract}
\maketitle
\tableofcontents
\section{Introduction}
One of the most fundamental theorems in Riemannian geometry is the Myers theorem~\cite{SBM}, which states that if an $n$-dimensional  complete Riemannian manifold $(M,g)$ satisfies $\Ric_g \geq (n-1)K$ with $K>0$ and $n\geq2,$ then $M$ is compact and its diameter $\diam(M)$ is at most $\pi/\sqrt{K}$. Moreover, its fundamental group $\pi_1 (M)$ is a finite group.
Here, we denote by $\Ric_g$ the Ricci curvature of $(M,g)$. For $\kappa \in C^\infty(M)$, $\Ric_g \geq \kappa$ means that $\Ric_g(v,v) \geq \kappa(p) g(v,v)$ holds for $v \in T_pM$. 
This theorem has been generalized by several ways. As one of these attempts, we introduce the theorems by Sprouse~\cite{SC} below. 

In this article, we always assume that $(M,g)$ is a complete Riemannian manifold and its dimension is at least $2$. We denote by $d$ and $\vol_g$ the Riemannian distance function and the volume measure on $(M,g)$, respectively. We write $a_+= \max\{a, 0\}$ for $a\in\R$. Set $\N_{\geq2}=\{\,n\in \N \mid n\geq2\,\}$. 
For $p\in M$ and $R>0$, we denote
\begin{align*}
U_p M &= \{\,v \in T_p M \mid g(v,v)=1\,\},\\
\Ric_{g-}(p) &= \inf \{\,\Ric_g(v,v) \mid v \in U_pM\,\},\\ 
B(p,R) &=\{\,q\in M\mid d(p,q)<R\,\}.
\end{align*}
\begin{thm}[\cite{SC}*{Theorem~1.1}]\label{Sp-cpt}
Let $n\in \N_{\geq2}$ and $\eta > 0$. There exists a positive constant  $\delta(n,\eta)$ such that if an $n$-dimensional compact Riemannian manifold $(M,g)$ satisfies $\Ric_g \geq 0$ and
\begin{equation*}
\frac{1}{\vol_g(M)} \int_M \big((n-1) - \Ric_{g-}\big)_+\,d\vol_g < \delta(n,\eta) 
\end{equation*} 
then $\diam(M) < \pi + \eta$.
\end{thm}

\begin{thm}[\cite{SC}*{Theorem~1.2}]\label{Sp-cmp}
Let $n\in \N_{\geq2}, K \leq 0$ and $R,\eta > 0$. There exists a positive constant $\delta (n,K,R,\eta)$ such that if an $n$-dimensional complete Riemannian manifold $(M,g)$ satisfies $\Ric_g\geq (n-1)K$ and
\begin{equation*}
\sup_{p\in M}\frac{1}{\vol_g(B(p,R))} \int_{B(p,R)} \big((n-1) - \Ric_{g-}\big)_+\, d\vol_g < \delta(n,K,R,\eta)
\end{equation*} 
then $M$ is compact and $\diam(M) < \pi + \eta$.  
\end{thm}

Hwang--Lee~\cite{HL} extended Theorems \ref{Sp-cpt} and \ref{Sp-cmp} to a weighted manifold having a bounded $\infty$-Ricci curvature. A \emph{weighted manifold} is a triple $(M,g,\mu_f )$, where $(M,g)$ is an $n$-dimensional Riemannian manifold, endowed with a measure $\mu_f = e^{-f}d\vol_g$ having smooth positive density with respect to the volume measure $\vol_g$. We call $f \in C^{\infty}(M)$ the \emph{weight function} of $\mu_f$. 
For given $N \in (-\infty,\infty]$, the \emph{$N$-Ricci curvature} is defined for $v\in TM$ by

\begin{align*}
\Ric_N (v,v)= 
 \begin{cases}
 \Ric_g(v,v) + \hess_g f(v,v) -\dfrac{1}{N-n}g(\nabla f, v)^2 & \text{if } N \in \R \setminus \{n\}, \\
 \Ric_g (v,v) + \hess_g f(v,v) & \text{if } N =\infty,
 \end{cases}
\end{align*}
and for $N=n$, 
\begin{align*}
\Ric_n(v,v) =
 \begin{cases}
\Ric_g(v,v) &\text{if } f \text{ is a constant function}, \\
-\infty & \text{otherwise}.
 \end{cases}
\end{align*}
We call the parameter $N$ the \emph{effective dimension}. We refer to a \emph{weighted Ricci curvature} as a generic term for the $N$-Ricci curvature with $N\in (-\infty,\infty]$. As similar as $\Ric_-$, we define $\Ric_{N-} : M \to \R$ by
\begin{equation*}
\Ric_{N-}(p)=\inf \{\,\Ric_N(v,v) \mid v \in U_pM\,\}. 
\end{equation*}

\begin{thm}[\cite{HL}*{Theorem~1.4}]
\label{HL-cpt}
Let $n\in \N_{\geq2}, k \geq 0$ and $\eta > 0$. 
There exists a positive constant $\delta (n+4k,\eta)$ such that if an $n$-dimensional compact weighted manifold $(M,g,\mu_f )$ satisfies $\Ric_\infty \geq 0$, $|f| \leq k$ and 
\begin{equation*}
\frac{1}{\mu_f (M)} \int_M \big((n-1) - \Ric_{\infty-} \big)_+ \, d\mu_f < \delta(n+4k,\eta)
\end{equation*} 
then
\begin{equation*}
\diam(M) < \lf(\pi + \frac{\eta}{2}\r)\sqrt{1+\frac{8k}{(n-1)\pi}} + \frac{\eta}{2}.
\end{equation*}
\end{thm}
\begin{thm}[\cite{HL}*{Theorem~1.5}]
\label{HL-cmp}
Let $n\in \N_{\geq2}, k\geq0, K \leq0$ and $R,\eta > 0$. There exists a positive constant $\delta (n+4k,K,R,\eta)$ such that if an $n$-dimensional complete weighted manifold $(M,g,\mu_f )$ satisfies $|f| \leq k$, $\Ric_\infty \geq (n-1)K$ and
\begin{equation*}
\sup_{p\in M}\frac{1}{\mu_f \big(B(p,R)\big)} \int_{B(p,R)}\big((n-1) - \Ric_{\infty-} \big)_+\, d\mu_f < \delta(n+4k,K,R,\eta) 
\end{equation*} 
then $M$ is compact and
\begin{equation*}
\diam(M) < \lf(\pi + \frac{\eta}{2}\r)\sqrt{1+\frac{8k}{(n-1)\pi}} + \frac{\eta}{2}.
\end{equation*}
\end{thm}

The aim of this article is to extend Theorems~\ref{Sp-cpt}--\ref{HL-cmp} for a weighted manifold under the assumption that the weighted curvature is bounded below in terms of its weight function. The key of the proof is the notion of the following $\vep$-range introduced by Lu--Minguzzi--Ohta~\cite{LMO}.
\begin{defi}
For $n \in \mathbb{N}_{\geq 2}$ and $N \in (-\infty,1] \cup [n,\infty]$, 
we say that \emph{$\vep \in \mathbb{R}$  is in the $\vep(n,N)$-range} if $\varepsilon \in \R$ satisfies the following conditions:
\begin{equation*}
 \varepsilon=0 \,\text{ for } N=1, \qquad
 \vert\varepsilon\vert < \sqrt{\frac{N-1}{N-n}} \,\text{ for } N \neq 1,n, \qquad
 \varepsilon \in \R \,\text{ for } N=n.
\end{equation*}
For $n,N$ as above and $\vep$ in the $\vep(n,N)$-range, we define the constant $c =c(n,N,\varepsilon)$ by
\begin{equation*}
c =\frac{1}{n-1}\lf(1-\varepsilon^2\frac{N-n}{N-1}\r) >0
\end{equation*}
for $N \neq 1$.
If $\varepsilon=0$, then one can take $N \to 1$ and set $c(1,0)=1/(n-1)$.
\end{defi}
Notice that $|\vep|=1$ happens only if $N \in [n,\infty)$. 
For $N \in (-\infty,1] \cup [n,\infty]$ and $K\in \R$, Lu--Minguzzi--Ohta~\cite{LMO2} established a curvature bound $\Ric_N \geq Ke^{4(\varepsilon -1)f/(n-1)}$ on a weighted manifold $(M,g,\mu_f)$ by using the $\vep(n,N)$-range. This curvature bound is a generalization of a different kind of curvature bounds $\Ric_1 \geq Ke^{4f/(1-n)}$ introduced by Wylie--Yeroshkin~\cite{WY} and for $N \in (-\infty,1]$, $\Ric_N \geq Ke^{4f/(N-n)}$ introduced by Kuwae--Li~\cite{KL}. They presented several comparison theorems with each curvature bounds. Furthermore, Kuwae--Sakurai~\cite{KS} also provided several comparison theorems with the $\vep(n,N)$-range. 

To state our theorems, we prepare the following condition.
\begin{defi}
Let $n\in \N_{\geq2}$, $N \in (-\infty,1] \cup [n,\infty]$, $\varepsilon$ 
in the $\vep(n,N)$-range, $K \leq 0$ and $b \geq a >0$.
\begin{enumerate}
\item 
We say that  an $n$-dimensional weighted manifold $(M,g,\mu_f )$ satisfies a \emph{$(N,K,\vep,a,b)$-condition} if one has
\begin{align*}
\Ric_N \geq Ke^{\frac{4(\varepsilon -1)}{n-1}f},
\qquad
a \leq e^{\frac{2(1-\varepsilon)}{n-1}f} \leq b\quad \text{on } M.
\end{align*}  
\item
We define the constants $\widetilde{D}=\widetilde{D}(n,N,\vep,a,b)$ by
\begin{align*}
&\widetilde{D}=
\begin{cases}
\sqrt{\dfrac{N-1}{n-1}} & \text{if } \vep=1, \\[-0.1cm]\\
\sqrt{1+\dfrac{2}{|1-\vep|\pi}\log\dfrac{b}{a}} & \text{if } \vep\neq1 \text{ and } N\in[n, \infty], \\[-0.1cm]\\ 
\sqrt{\lf(\dfrac{b}{a}\r)^{\lmd_0} 
	+ \dfrac{2}{(1-\vep)\pi\lmd_0}\lf(\lf(\dfrac{b}{a}\r)^{\lmd_0}-1\r)} & \text{otherwise}.
 \end{cases}
\end{align*}
Moreover, for $N\in (-\infty,1]$, we define $\lmd_0=\lmd_0(n,N,\vep)$ by
\begin{equation*}
\lmd_0=\frac{1}{1-\vep}\lf(1-\sqrt{\frac{N-1}{N-n}}\r).
\end{equation*}
\end{enumerate}
\end{defi}

For a weighted manifold satisfying the $(N,K,\vep,a,b)$-condition, we estimate its diameter quantitatively.
\begin{thm}\label{Cpt-Main} 
Let $n\in \N_{\geq2}, N \in (-\infty,1] \cup [n,\infty]$, $\varepsilon$ in the $\vep(n,N)$-range, $b\geq a >0$ and $H,\eta>0$. 
There exists a positive constant $\delta_1(n,N,\vep,a,b,H,\eta)$ such that if an $n$-dimensional compact weighted manifold $(M,g,\mu_f )$ satisfies the $(N,0,\vep,a,b)$-condition and
\begin{equation*}
\frac{1}{\mu_f (M)} \int_M \big((n-1)H - \Ric_{N-}\big)_+ \,d\mu_f \leq \delta_1(n,N,\vep,a,b,H,\eta)
\end{equation*} 
then 
\begin{equation*}
\diam(M) \leq \frac{\pi+\eta}{\sqrt{H}}\widetilde{D}(n,N,\vep,a,b).
\end{equation*}
\end{thm}

\begin{thm}\label{Cmp}
Let $n\in \N_{\geq2}, N \in (-\infty,1] \cup [n,\infty]$ and $\varepsilon$ in the $\vep(n,N)$-range. Take real numbers $K, a, b, H, R$ and $\eta$ such that $K\leq 0$ and $a,b,H,R,\eta>0$ with $a\leq b$. There exists a positive constant $\delta(n,N,K,\vep,a,b,H,R,\eta)$ such that if an $n$-dimensional complete weighted manifold $(M,g,\mu_f )$ satisfies the $(N,K,\vep,a,b)$-condition and
\begin{equation*}
\sup_{p\in M} \frac{1}{\mu_f \big(B(p,R)\big)} \int_{B(p,R)} \big((n-1)H - \Ric_{N-}\big)_+ \,d\mu_f \leq \delta(n,N,K,\vep,a,b,H,R,\eta)
\end{equation*} 
then $M$ is compact and
\begin{equation*}
\diam(M) \leq \frac{\pi+\eta}{\sqrt{H}}\widetilde{D}(n,N,\vep,a,b).
\end{equation*}
\end{thm}
\begin{rem}
Assume that an $n$-dimensional complete weighted manifold $(M,g,\vol_g)$ satisfies $\Ric_g\geq (n-1)H$ for $H>0$. Then $(M,g,\vol_g)$ satisfies the $(n,0,1,1,1)$-condition and
\begin{equation*}
\sup_{p\in M} \frac{1}{\vol_g \big(B(p,R)\big)} \int_{B(p,R)} \big((n-1)H - \Ric_{n-}\big)_+ \,d\vol_g =0
\end{equation*}
holds for each $R>0$.  Furthermore, we find $\widetilde{D}(n,n,1,1,1)=1$. By Theorem~\ref{Cmp}, we see that $\diam(M)\leq (\pi+\eta)/\sqrt{H}$ for any $\eta>0$. Taking $\eta\to 0$, we obtain $\diam(M) \leq \pi/\sqrt{H}$. This means that Theorem~\ref{Cmp} recovers the Myers theorem.
\end{rem}
This article is organized as follows: In Section~\ref{PS}, we recall a Bishop-type inequality, the volume comparison theorem on a weighted manifold given in~\cite{LMO2}. Then we extend the segment inequality of Cheeger--Colding~\cite{CC} to a weighted manifold under the assumption that the weighted Ricci curvature is bounded below in terms of its weight function. In Sections~\ref{M-compact} and~\ref{M-complete}, we prove Theorems~\ref{Cpt-Main} and~\ref{Cmp}, respectively. In Sections~\ref{FGB}, we analyze the fundamental group of a weighted manifold satisfying the $(N,K,\vep,a,b)$-condition.
%

\begin{Acknowledgements}
The author would like to express his deepest thanks to his supervisor, Asuka Takatsu, as well as to Manabu Akaho, Shin-ichi Ohta, Takashi Sakai, Homare Tadano who contributed their support. 
\end{Acknowledgements}
%
%
\section{Preliminaries and segment inequality  with $\vep$-range}\label{PS}
%
%
We recall some results by Lu--Minguzzi--Ohta~\cite{LMO2}. Although they discussed a weighted Finsler manifold, throughout this article, we treat a weighted manifold.
%
\subsection{Preliminaries}
%
Lu--Minguzzi--Ohta~\cite{LMO2} introduced the $\vep$-range and  proved a volume comparison theorem for a weighted manifold $(M,g,\mu_f)$.

For $p\in M$, let $\gamma : [0,l) \to M$ be a unit speed geodesic such that $\gamma(0)=p$. Set $f_\gamma (t)=f\big(\gamma(t)\big)$. We denote by 
\begin{equation*}
\mu_f =e^{-f}d\vol_g= \A^f_{\gamma}(t)dt d\theta_{n-1} =e^{-f_\gamma(t)}\A_{\gamma}(t)dt d\theta_{n-1}
\end{equation*}
the weighted measure in the geodesic polar coordinates along geodesics $\gamma$, where $d\theta_{n-1}$ is the volume measure of the unit sphere in $T_pM$. 

First we introduce a Bishop-type inequality.

\begin{prop}[\cite{LMO2}*{Theorem~3.5}]
\label{Bishop}
Let $(M,g,\mu_f )$ be an $n$-dimensional complete weighted manifold, $N \in (-\infty,1] \cup [n,\infty]$, $\varepsilon$ 
in the $\vep(n,N)$-range and $c=c(n,N,\vep)$. 
For a unit speed geodesic $\gamma : [0,l) \to M$, set 
\begin{equation*}
 h(t) = e^{-cf_\gamma(t)}\A_{\gamma}(t)^c, \qquad
 h_1(\tau)=h ( \varphi_\gamma^{-1}(\tau) )
 \end{equation*}
for $t \in [0,l)$ and $\tau \in [0,\varphi_\gamma(l))$, where
\begin{equation*}\label{eq:phi_e}
\varphi_\gamma(t) =\int_0^t e^{\frac{2(\varepsilon -1)}{n-1}f_\gamma(s)} \, ds.
\end{equation*}
Then, for all $\tau \in (0,\varphi_\gamma(l))$, 
\begin{equation*}\label{eq:F-Bish}
h_1''(\tau) \leq -ch_1(\tau)\Ric_N \big( (\gamma\circ\varphi_\gamma^{-1} )' (\tau) \big).
\end{equation*}
\end{prop}

We define the comparison function $\sn_\kappa$ by
\begin{equation*}
\sn_{\kappa}(t) = \begin{cases}
 \dfrac{1}{\sqrt{\kappa}} \sin(\sqrt{\kappa}t) & \kappa>0, \\
 t &\kappa=0,  \\
 \dfrac{1}{\sqrt{-\kappa}} \sinh(\sqrt{-\kappa}t) & \kappa<0,
 \end{cases}
\end{equation*}
where $t \in [0,\pi/\sqrt{\kappa}]$ for $\kappa>0$ and $t \in \R$ for $\kappa\leq 0$. Notice that $\sn_\kappa$ satisfies 
\begin{equation*}
 \begin{cases}
	\sn_\kappa''(t) + \kappa\sn_\kappa(t)=0,\\
	\sn_\kappa(0)=0,\\
	\sn'_\kappa(0)=1.
 \end{cases}
\end{equation*}

Let us recall a volume comparison theorem.
\begin{prop}[\cite{LMO2}*{Theorem~3.11}]\label{BCT}
Let $(M,g,\mu_f )$ be an $n$-dimensional complete weighted manifold satisfying the $(N,K,\vep,a,b)$-condition. Then
\begin{equation*}
\dfrac{\mu_f \big(B(p,R)\big)}{\mu_f \big(B(p,r)\big)} 
\leq \dfrac{b}{a} \cdot \dfrac{\displaystyle\int_0^{R/a} \sn_{cK}(\tau)^{1/c}\,d\tau}{\displaystyle\int_0^{r/b} \sn_{cK}(\tau)^{1/c}\,d\tau}
\end{equation*} 
holds for all $p \in M$ and $0 < r \leq R$.  
\end{prop}
It should be mentioned that Theorem~\ref{BCT} was originaly proved for all $K\in \R$ in~\cite{LMO2}*{Theorem~3.11}.
\begin{rem}
If $K=0$, we have $\sn_0(\tau) =\tau$, which gives
\begin{align*}
\int_0^{r} \sn(\tau)^{\frac{1}{c}}\, d\tau
= \frac{c}{c+1}r^{\frac{1}{c}+1}
\end{align*}
for $r>0$. In the setting of Theorem~\ref{BCT} for $K=0$, we find that
\begin{equation*}
\dfrac{\mu_f \big(B(p,R)\big)}{\mu_f \big(B(p,r)\big)} 
\leq \dfrac{b}{a} \cdot \lf(\frac{bR}{ar}\r)^{\frac{1}{c}+1}.
\end{equation*}  
\end{rem}
\subsection{Segment inequality with $\vep$-range}
%
The segment inequality proved by Cheeger--Colding~\cite{CC}*{Theorem~2.11} plays an important role in the proof  of Theorems~\ref{Sp-cpt} and~\ref{Sp-cmp}. We extend this theorem as follows. 
\begin{thm}\label{Seg-Ran}
Let $(M,g,\mu_f )$ be an $n$-dimensional complete weighted manifold, $N\in (-\infty,1] \cup [n,\infty]$, $\varepsilon$ in the $\vep(n,N)$-range and $K \in \R$. Assume that   
\begin{equation*}
\Ric_N \geq Ke^{\frac{4(\varepsilon -1)}{n-1}f}
\end{equation*} 
holds.
For $i=1,2$, let $A_i$ be bounded open subsets of $M$ and W be an open subset of $M$ such that, for each two points $y_i \in A_i$, any unit minimal geodesics $\gamma_{y_1,y_2}$ from $y_1$ to $y_2$ is containd in W.
Then for any non-negative integrable function $F$ on $W$, we have
\begin{align*}
&\quad \int_{A_1\times A_2}\lf(\int_0^{d(y_1,y_2)} F(\gamma_{y_1, y_2})\,ds \r) d\mu_{f\times f} \\
&\leq C(n,c,K)\lf[\mu_f (A_2)\,\diam(A_1)+\mu_f (A_1)\,\diam(A_2)\r]\int_W F\, d\mu_f ,
\end{align*} 
where $\mu_{f\times f}$ is the product measure on $M\times M$ induced by $\mu_f $ and
\begin{equation*}
C(n,c,K)
= \sup_{y_1 \in A_1, y_2\in A_2}
\lf(\sup_{0 < \frac{s}{2} \leq u < s\leq d(y_1,y_2)}
\frac{\sn_{cK} \big(\varphi_\gamma(s)\big)^{1/c}}{\sn_{cK} \big(\varphi_\gamma(u)\big)^{1/c}}\r).
\end{equation*}
\end{thm}
\begin{pf}
Set 
\begin{equation*}
B= \{(y_1,y_2) \in A_1 \times A_2 \mid \mathrm{there\ exists\ a\ unique\ minimal\ geodesic\ from}\ y_1\ \mathrm{to}\ y_2\}.
\end{equation*}
Then $\mu_{f\times f}(B) =\mu_{f\times f}(A_1 \times A_2)$ holds since the measure of a cut locus of each $y_i \in A_i$ with respect to $\mu_f$ is zero for $i=1,2$.

Define maps $E_1,E_2,E : A_1\times A_2 \to \R$ by 
\begin{align*}
E_1(y_1,y_2) &= \int_{d(y_1,y_2)/2}^{d(y_1,y_2)}  F\big(\gamma_{y_1,y_2}(u)\big)\,du , \\
E_2(y_1,y_2) &= \int_0^{d(y_1,y_2)/2} F\big(\gamma_{y_1,y_2}(u)\big)\,du ,\\
E(y_1,y_2) &= E_1(y_1,y_2) +E_2(y_1,y_2) = \int_0^{d(y_1,y_2)} F\big(\gamma_{y_1,y_2}(u)\big)\,du .
\end{align*}

Fix $y_i \in A_i$ and $v_i \in U_{y_i}M$ for $i=1,2$. Set 
\begin{equation*}
I(y_i,v_i) =\{\,s>0\mid \exp_{y_i}(s v_i) \in A_{i+1}, \,d\big(y_i,\exp_{y_i}(s v_i)\big) =s\,\},
\end{equation*}
where we put $A_{3}=A_1$. We denote by $|I(y_1,v_1)|$ the $1$-dimensional Lebesgue measure of $I(y_1,v_1)$.  Since $\{\,\exp_{y_i}(s v_i) \mid s\in I(y_i,v_i)\,\}$ is contained only in $A_{i+1}$, we have
\begin{equation*}
|I(y_i,v_i)| \leq \diam(A_{i+1})
\end{equation*}
for $i=1,2$. 

Fix $s \in I(y_1,v_1)$. For $\gamma(u)=\exp_{y_1}(u v_1)$ on $u\in [0,s)$, by Proposition~\ref{Bishop}, we find that $\mathcal{A}^f_{\gamma}(u)/\big(\sn_{cK} (\varphi_\gamma(u))\big)^{1/c}$ is non-increasing, where $\varphi_\gamma(u)<\pi/\sqrt{cK}$ holds if $K>0$ (see~\cite{LMO2}*{Proof of Theorem~3.6}). Then we find that
\begin{equation*}
\frac{\mathcal{A}^f_{\gamma}(s)}{\mathcal{A}^f_{\gamma}(u)}
\leq \frac{\sn_{cK} \big(\varphi_\gamma(s)\big)^{1/c}}{\sn_{cK} \big(\varphi_\gamma(u)\big)^{1/c}}
\leq
C(n,c,K)
\end{equation*}
for $0<s/2 \leq u< s\leq T(v_1)$, where $T(v_1)$ is the supremum of $s$ such that $ s\in I(y_1, v_1)$. This yields
\begin{align*}
E_1(y_1,\gamma(s))\mathcal{A}^f_{\gamma}(s)
&= \mathcal{A}^f_{\gamma}(s) \int^s_{s/2} F\big(\gamma(u)\big) \,du \\
&\leq C(n,c,K) \int^s_{s/2} F\big(\gamma(u)\big) \mathcal{A}^f_{\gamma}(u)\,du\\
&\leq C(n,c,K) \int^{T(v_1)}_0 F\big(\gamma(u)\big) \mathcal{A}^f_{\gamma}(u)\,du.
\end{align*}
Integrating this inequality on $I(y_1,v_1)$ implies that
\begin{align*}
\int_{I(y_1,v_1)}E_1(y_1,\gamma(s))\mathcal{A}^f_{\gamma}(s) \,ds
&\leq C(n,c,K) \int_{I(y_1,v_1)} \,ds \int_0^{T(v_1)} F\big(\gamma(u)\big) \mathcal{A}^f_{\gamma}(u)\,du \\
&\leq C(n,c,K)\, \diam(A_2) \int_0^{T(v_1)} F\big(\gamma(u)\big) \mathcal{A}^f_{\gamma}(u)\,du. 
\end{align*}
By integrating this inequality over the unit sphere in $T_{y_1}M$, we see that
\begin{align*}
\int_{A_2}E_1(y_1,y_2)d\mu_f 
&= \int_{U_{y_1}M}\,d\theta_{n-1} \int_{I(y_1,v_1)}E_1(y_1,\gamma(s))\mathcal{A}^f_{\gamma}(s) \,ds \\
&\leq C(n,c,K)\, \diam(A_2)  \int_{U_{y_1}M}\,d\theta_{n-1} \int_0^{T(v_1)} F\big(\gamma(u)\big) \mathcal{A}^f_{\gamma}(u)\,du \\
&\leq C(n,c,K)\, \diam(A_2)  \int_W F \,d\mu_f ,
\end{align*}
where, in the third inequality, we use the fact that any unit minimal geodesics from $y_1 \in A_1$ to $y_2\in A_2$ is containd in $W$. Therefore we have
\begin{equation}\label{PAA}
\int_{A_1\times A_2}E_1(y_1,y_2)d\mu_{f\times f}
\leq C(n,c,K)\, \diam(A_2)\, \mu_f (A_1)\int_W F \,d\mu_f .
\end{equation}
On the other hand, substituting $t=d(y_1,y_2)-u$ for $E_2(y_1,y_2)$ yields
\begin{align}\label{chg}
\begin{split}
E_2(y_1,y_2)
&=\int_0^{d(y_1,y_2)/2} F\big(\gamma_{y_1,y_2}(u)\big)\,du \\
&=\int_{d(y_2,y_1)/2}^{d(y_2,y_1)} F\big(\sigma_{y_2,y_1}(t)\big)\,dt \\  
&=E_1(y_2,y_1),
\end{split}
\end{align}
where $\sigma_{y_2,y_1}(t)=\gamma_{y_1,y_2}\big(d(y_1,y_2)-t\big)$ on $[d(y_2,y_1)/2, d(y_2,y_1)]$. Interchanging the roles of $A_1$ and $A_2$ in \eqref{PAA}, it turns out that, by \eqref{chg}, 
\begin{align*}
\int_{A_1\times A_2}E_2(y_1,y_2)\,d\mu_{f\times f}
&=\int_{A_1\times A_2}E_1(y_2,y_1)\,d\mu_{f\times f} \\
&\leq C(n,c,K)\, \diam(A_1)\, \mu_f (A_2)\int_W F \,d\mu_f . 
\end{align*}
Thus we obtain
\begin{align*}
&\quad \int_{A_1\times A_2}E(y_1,y_2)\,d\mu_{f\times f}\\
&= \int_{A_1\times A_2}E_1(y_1,y_2)\,d\mu_{f\times f} + \int_{A_1\times A_2}E_2(y_1,y_2)\,d\mu_{f\times f} \\
&\leq C(n,c,K) \big[\diam(A_2)\, \mu_f (A_1)+ \diam(A_1)\, \mu_f (A_2)\big]\int_W F \,d\mu_f .  
\end{align*}
This completes the proof of the theorem.
\end{pf}
\begin{rem}
Set $S = \sup_{y_1 \in A_1, y_2\in A_2}d(y_1,y_2)$. If we take $N=n$, $f \equiv 0$, $\varepsilon =0$ and $K=(n-1)H$ with $H>0$, then we have
\begin{equation*}
C(n,c,K)
= \sup_{0 < \frac{s}{2} \leq u < s\leq S}
\lf(\frac{\sn_{H}(s)}{\sn_{H} (u)}\r)^{n-1}
\end{equation*}
and the estimate in Theorem~\ref{Seg-Ran} coincides with that of~\cite{CC}*{Theorem~2.11}. Whereas, for $N= \infty$, since we impose the different condition compared with~\cite{JM}*{Proposition 2.3}, the estimate in Theorem~\ref{Seg-Ran} differs from that of~\cite{JM}*{Proposition 2.3}   
\end{rem}
Finally, we provide the following lemma towards the proof of Theorems~\ref{Cpt-Main} and~\ref{Cmp}.
%
\begin{lem}\label{Inc-sn}
Let $(M,g,\mu_f )$ be an $n$-dimensional complete weighted manifold satisfying the $(N,K,\vep,a,b)$-condition. Take $A_1$, $A_2$ and $C(n,c,K)$ as in \rm{Theorem~\ref{Seg-Ran}}. Then 
\begin{equation*}
C(n,c,K)\leq
\frac{\sn_{cK} (S/a)^{1/c}}{\sn_{cK} (S/2b)^{1/c}}
\end{equation*}
holds, where $S=\sup_{y_1 \in A_1, y_2\in A_2} d(y_1,y_2)$. 
\end{lem}
\begin{pf}
The boundedness of the weight function $f$ yields $s/b \leq \varphi_\gamma(s) \leq s/a$ for any unit speed minimal geodesic $\gamma$ from $y_1 \in A_1$ to $y_2\in A_2$. Then for $K\leq0$, the monotonicity of $\sn_{cK}(\tau)$ in $\tau >0$ gives
\begin{equation}\label{clm}
C(n,c,K)
\leq 
\sup_{0 < s\leq S} \frac{\sn_{cK} (s/a)^{1/c}}{\sn_{cK} (s/2b)^{1/c}}.
\end{equation}
For $K=0$, $\sn_{0}(s)=s$ implies
\begin{equation*}
C(n,c,K)
\leq \sup_{0 < s\leq S} \frac{\sn_{0} (s/a)^{1/c}}{\sn_{0} (s/2b)^{1/c}}
= \lf(\frac{2b}{a}\r)^{\frac{1}{c}}.
\end{equation*}

When $K<0$, we will see the following claim:
\begin{claim}\label{sn-clm}Let $F : \R \to \R$ be
\begin{equation*}
F(s)=\frac{\sinh(As)}{\sinh(Bs)}
\end{equation*}
for $A>B>0$. Then $F(s)$ is strictly increasing on $(0,\infty)$.
\end{claim}
\begin{cpf}
Differentiating $F(s)$ gives
\begin{equation*}
F'(s) 
= \frac{1}{\sinh^2(Bs)}
\lf(A\cosh(As)\sinh(Bs)- B\sinh(As)\cosh(Bs)\r).
\end{equation*}
We use the following hyperbolic function identities:
\begin{align*}
2\cosh(As)\sinh(Bs)
= \sinh\lf(A+B\r)s + \sinh\lf(-A+B\r)s,\\
2\sinh(As)\cosh(Bs)
= \sinh\lf(A+B\r)s - \sinh\lf(-A+B\r)s.
\end{align*}
These yield
\begin{equation*}
F'(s) 
= \frac{1}{2\sinh^2(Bs)}\lf\{
\lf(A-B\r)\sinh\lf(A+B\r)s- \lf(A+B\r)\sinh\lf(A-B\r)s\r\} .
\end{equation*}
Setting
\begin{equation*}
F_1(s)=
\lf(A-B\r)\sinh\lf(A+B\r)s- \lf(A+B\r)\sinh\lf(A-B\r)s,
\end{equation*}
we have $F_1(0)=0$ and 
\begin{equation*}
F'_1(s)=
AB\lf\{
\cosh\lf(A+B\r)s- \cosh\lf(A-B\r)s\r\}> 0
\end{equation*}
on $(0,\infty)$. Therefore we get $F_1(s) > 0$ on $(0,\infty)$. Hence we obtain $F'(s)>0$ on $(0,\infty)$. This completes the proof of the claim.
\end{cpf}
When $K< 0$, applying this claim to \eqref{clm} gives
\begin{equation*}
C(n,c,K)
\leq \frac{\sn_{cK} (S/a)^{1/c}}{\sn_{cK} (S/2b)^{1/c}}.
\end{equation*}
\end{pf}
%
\section{Proof of Theorem~\ref{Cpt-Main}}\label{M-compact}
%
We extend Theorem \ref{Sp-cpt} to the case that the weighted manifold has the non-negative weighted curvature.  
\begin{apf}{Theorem~\ref{Cpt-Main}}
Set $D=\diam(M)$ and take $p_1,p_2 \in M$ such that $D=d(p_1,p_2)$.
We put $W=M$ and $A_i=B(p_i,r)$ for $i=1,2$, where $r>0$ is later determined.
Then Theorem~\ref{Seg-Ran} gives 
\begin{align}
\begin{split}
&\quad \int_{A_1\times A_2}\lf(\inf_{(y_1,y_2) \in \overline{A}_1 \times \overline{A}_2} \int_0^{d(y_1,y_2)} \big((n-1)H - \Ric_{N-}\big)_+(\gamma_{y_1, y_2})\,ds \r) d\mu_{f\times f} \\
&\leq \int_{A_1\times A_2}\lf(\int_0^{d(y_1,y_2)} \big((n-1)H - \Ric_{N-}\big)_+(\gamma_{y_1, y_2})\,ds \r) d\mu_{f\times f} \\
&\leq 2rC(n,c,K)\lf[\mu_f (A_1)+\mu_f (A_2)\r]\int_M \big((n-1)H - \Ric_{N-}\big)_+\, d\mu_f \label{Segment}. 
\end{split}
\end{align} 
We observe from Theorem~\ref{BCT} with $K=0$ that
\begin{equation*}
\dfrac{\mu_f (M)}{\mu_f (A_i)}
\leq \frac{b}{a}\lf(\frac{bD}{ar}\r)^{\frac{1}{c}+1}
\end{equation*}
for $i=1,2$ and from Lemma~\ref{Inc-sn} that
\begin{equation*}
C(n,c,K)\leq \lf(\frac{2b}{a}\r)^{\frac{1}{c}}.
\end{equation*} 
Dividing \eqref{Segment} by $\mu_f (A_1)\mu_f (A_2)$ yields 
\begin{align}
\begin{split}\label{L-com}
&\quad \inf_{(y_1,y_2) \in \overline{A}_1 \times \overline{A}_2}
\int_0^{d(y_1,y_2)} \big((n-1)H - \Ric_{N-}\big)_+(\gamma_{y_1,y_2})\,ds \\
&\leq 2rC(n,c,K)\lf(\frac{1}{\mu_f (A_1)}+\frac{1}{\mu_f (A_2)}\r)\int_M \big((n-1)H - \Ric_{N-}\big)_+\, d\mu_f \\
&\leq 2r\,\lf(\frac{2b^2D}{a^2r}\r)^{\frac{1}{c}+1}\frac{1}{\mu_f (M)}\int_M \big((n-1)H - \Ric_{N-}\big)_+\, d\mu_f . 
\end{split}
\end{align} 

There exists a unit speed minimal geodesic $\gamma$ from $y_1 \in \overline{A}_1$ to $y_2 \in \overline{A}_2$ that
attains the infimum of \eqref{L-com}. 
Let $L=d(y_1,y_2)$ and $\{E_1,\ldots,E_n = \dot{\gamma}\}$ be a parallel orthonormal frame along $\gamma$. For a smooth function $\alpha \in C^\infty([0,L])$ such that $\alpha(0)=\alpha(L)=0$, we set $Y_i(t) = \alpha(t)E_i(t)$, $i=1,\ldots,n-1$.
We denote by $L_i(s)$ the length functional of a fixed-endpoint variation of a curves $c(s,t) :(-\epsilon,\epsilon)\times[0,L] \to M$ such that 
\begin{equation*}
c(0,t)=\gamma(t), \qquad \lf.\frac{\partial}{\partial s} c(s,t)\r|_{s=0}=Y_i(t). 
\end{equation*}
Then the second variation formula for $L_i(s)$ (see~\cite{ST}*{Chapter III Theorem 2.5}) provides
\begin{align}\label{com-totyuu}
\lf.\sum_{i=1}^{n-1}\frac{d^2L_i}{ds^2}\r|_{s=0}
&= \sum_{i=1}^{n-1}\int_0^L \lf\{g(\nabla_{\dot{\gamma}}Y_i, \nabla_{\dot{\gamma}}Y_i) - R_g\big(Y_i(t),\gamma'(t),\gamma'(t) ,Y_i(t)\big) \r\}\,dt \\
&= \int_0^L \lf\{(n-1)\alpha'(t)^2 - \alpha(t)^2\Ric_g(\gamma'(t),\gamma'(t)) \r\}\,dt \notag\\
&= \int_0^L \big[-(n-1)H\alpha(t)^2 + (n-1)\alpha'(t)^2 +\alpha(t)^2\hess_g\,f\big(\gamma'(t),\gamma'(t)\big)  \notag\\
&\quad -\frac{\alpha(t)^2}{N-n}f'_\gamma(t)^2
	 + \alpha(t)^2\lf\{(n-1)H - \Ric_N\big(\gamma'(t),\gamma'(t)\big) \r\}\big]\,dt ,\notag 
\end{align}
where $R_g$ is the Riemannian curvature tensor of $(M,g)$. Note that if $N=n$, the forth term on the right-hand side of \eqref{com-totyuu} always vanishes since we set $f'_\gamma(t)^2/(N-n)=0$. 

We take a parameter $\lmd$ satisfying
\begin{align}\label{lmd-asmp}
\begin{cases}
|(1-\vep)\lmd-1|\leq \sqrt{\dfrac{N-1}{N-n}} \qquad &\text{ if }N\in(-\infty,1]\cup(n,\infty] \text{ and }\vep\neq1,\\[-0.1cm]\\
\lmd \in \R \qquad &\text{ if either }N=n \text{ or } \vep=1.
\end{cases} 
\end{align}
We define $\Phi(\lmd)$, $\Psi(\lmd)$ by
\begin{equation*}
 \Phi(\lmd)=
 \begin{cases}
 b^\lmd \, &\text{if } \lmd>0, \\
 1\, &\text{if } \lmd=0, \\
 a^\lmd\, &\text{if } \lmd<0,
 \end{cases} \qquad
 \Psi(\lmd) =
  \begin{cases}
 a^\lmd \, &\text{if } \lmd>0, \\
 1\, &\text{if } \lmd=0, \\
 b^\lmd\, &\text{if } \lmd<0.
 \end{cases}
\end{equation*}
If we choose 
\begin{equation*}
\alpha(t) = e^{\frac{(1-\vep)\lmd}{n-1}f_\gamma(t)} \sin \lf(\frac{\pi t}{L}\r),
\end{equation*} 
then we have
\begin{equation*}
\alpha'(t) = \frac{(1-\vep)\lmd}{n-1} f'_\gamma(t) \alpha(t) +\frac{\pi t}{L} e^{\frac{(1-\vep)\lmd}{n-1}f_\gamma(t)}\cos \lf(\frac{\pi t}{L}\r).
\end{equation*}
We estimate the right-hand side of \eqref{com-totyuu}.
The first term constructs to
\begin{align}
-(n-1)H\int_0^L \alpha(t)^2 \,dt 
&\leq -(n-1)H \int_0^L  \Psi(\lmd)\sin^2 \lf(\frac{\pi t}{L}\r) \,dt \label{fst-term}\\
&\leq  -\frac{(n-1)HL}{2} \Psi(\lmd). \notag 
\end{align}
The second term is estimated as
\begin{align}
&\quad \int_0^L (n-1)\alpha'(t)^2\,dt \label{sec-term} \\
&= \frac{(n-1)\pi^2}{L^2}\int_0^L e^{\frac{2(1-\vep)\lmd}{n-1}f_\gamma(t)} \cos^2 \lf(\frac{\pi t}{L}\r) \,dt 
	+\frac{(1-\vep)\pi\lmd}{L}\int_0^L f'_\gamma(t)e^{\frac{2(1-\vep)\lmd}{n-1}f_\gamma(t)} \sin \lf(\frac{2\pi t}{L}\r) \,dt \notag\\
	&\quad+\frac{(1-\vep)^2\lmd^2}{n-1}\int_0^L\alpha(t)^2 f'_\gamma(t)^2 \,dt \notag\\
&= \frac{(n-1)\pi^2}{L^2}\int_0^L e^{\frac{2(1-\vep)\lmd}{n-1}f_\gamma(t)} \cos^2 \lf(\frac{\pi t}{L}\r) \,dt 
	-\frac{(n-1)\pi^2}{L^2}\int_0^L e^{\frac{2(1-\vep)\lmd}{n-1}f_\gamma(t)} \cos \lf(\frac{2\pi t}{L}\r) \,dt \notag\\
	&\quad+\frac{(1-\vep)^2\lmd^2}{n-1}\int_0^L\alpha(t)^2 f'_\gamma(t)^2 \,dt \notag\\
&= \frac{(n-1)\pi^2}{2L^2}\int_0^L e^{\frac{2(1-\vep)\lmd}{n-1}f_\gamma(t)}\lf(1-\cos \lf(\frac{\pi t}{L}\r)\r) \,dt
	+\frac{(1-\vep)^2\lmd^2}{n-1}\int_0^L\alpha(t)^2 f'_\gamma(t)^2 \,dt\notag \\
&\leq \frac{(n-1)\pi^2}{2L} \Phi(\lmd)
	+\frac{(1-\vep)^2\lmd^2}{n-1}\int_0^L\alpha(t)^2 f'_\gamma(t)^2 \,dt. \notag
\end{align}
We calculate the third term as
\begin{align}
\begin{split}
&\quad \int_0^L \alpha(t)^2f_\gamma''(t) \,dt \label{thd-term}\\
&= -\int_0^L 2\alpha(t)\alpha'(t)f'_\gamma(t)\,dt \\
&= -\frac{\pi}{L}\int_0^L f'_\gamma(t)e^{\frac{2(1-\vep)\lmd}{n-1}f_\gamma(t)} \sin \lf(\frac{2\pi t}{L}\r) \,dt
	-\frac{2(1-\vep)\lmd}{n-1}\int_0^L \alpha(t)^2 f'_\gamma(t)^2  \,dt .
\end{split}
\end{align}
The last term constructs to
\begin{align}\label{last-term}
\begin{split}
 &\quad \int_0^L \alpha(t)^2\lf\{(n-1)H - \Ric_{N-}\big(\gamma'(t),\gamma'(t)\big) \r\}\,dt \\
&\leq \int_0^L \alpha(t)^2 \big((n-1)H - \Ric_{N-} \big)_+\,dt \\
&\leq \Phi(\lmd)\int_0^L \big((n-1)H - \Ric_{N-} \big)_+ \,dt.
\end{split}
\end{align}
Combining \eqref{fst-term}, \eqref{sec-term}, \eqref{thd-term}, \eqref{last-term} gives
\begin{align}
\lf.\sum_{i=1}^{n-1}\frac{d^2L_i}{ds^2}\r|_{s=0} 
&\leq -\frac{(n-1)HL}{2} \Psi(\lmd) 
	+ \frac{(n-1)\pi^2}{2L} \Phi(\lmd)\\
	 &\quad+ \widetilde{D}_1(\lmd) + \Phi(\lmd)\int_0^L \big((n-1)H - \Ric_{N-} \big)_+ \,dt , \notag
\end{align}
where
\begin{align*}
\widetilde{D}_1(\lmd) 
&= -\frac{\pi}{L}\int_0^L f'_\gamma(t)e^{\frac{2(1-\vep)\lmd}{n-1}f_\gamma(t)} \sin \lf(\frac{2\pi t}{L}\r) \,dt \\
	&\quad + \lf(\frac{(1-\vep)^2\lmd^2}{(n-1)^2} -\frac{2(1-\vep)\lmd}{n-1} -\frac{1}{N-n}\r)\int_0^L \alpha(t)^2 f'_\gamma(t)^2  \,dt.
\end{align*}
First we consider the case of $\vep=1$. We see that
\begin{align*}
\widetilde{D}_1(\lmd)
&= -\frac{\pi}{L}\int_0^L f'_\gamma(t) e^{\frac{2(1-\vep)\lmd}{n-1}f_\gamma(t)}\sin \lf(\frac{2\pi t}{L}\r) \,dt - \frac{1}{N-n} \int_0^L \alpha(t)^2 f'_\gamma(t)^2  \,dt \\
&= (N-n)\frac{\pi^2}{L^2}\int_0^L e^{\frac{2(1-\vep)\lmd}{n-1}f_\gamma(t)}\cos^2 \lf(\frac{\pi t}{L}\r) \,dt \\
&\quad - (N-n)\int_0^Le^{\frac{2(1-\vep)\lmd}{n-1}f_\gamma(t)}\lf(\frac{\pi}{L}\cos\lf(\frac{\pi t}{L}\r) +\frac{f'_\gamma(t)}{N-n}\sin\lf(\frac{\pi t}{L}\r)\r)^2 \,dt \\
&\leq \frac{(N-n)\pi^2}{2L}\Phi(\lmd).
\end{align*}
If $\vep\neq 1$, then $\lmd\neq 0$ and it follows from \eqref{lmd-asmp} that
\begin{align*}
\widetilde{D}_1(\lmd) 
&= -\frac{\pi}{L}\int_0^L f'_\gamma(t)e^{\frac{2(1-\vep)\lmd}{n-1}f_\gamma(t)} \sin \lf(\frac{2\pi t}{L}\r) \,dt \\
	&\quad + \lf(\frac{(1-\vep)^2\lmd^2}{(n-1)^2} -\frac{2(1-\vep)\lmd}{n-1} -\frac{1}{N-n}\r)\int_0^L \alpha(t)^2 f'_\gamma(t)^2  \,dt \\
&\leq \frac{(n-1)\pi^2}{(1-\vep)L^2\lmd}\int_0^L e^{\frac{2(1-\vep)\lmd}{n-1}f_\gamma(t)} \cos \lf(\frac{2\pi t}{L}\r) \,dt \\
&\leq \frac{(n-1)\pi}{|1-\vep|L\lmd}\big(\Phi(\lmd) - \Psi(\lmd)\big).
\end{align*}
Moreover, if we define
\begin{align*}
\widetilde{D}(\lmd)=
\begin{cases}
\sqrt{\dfrac{N-1}{n-1}\lf(\dfrac{\Phi(\lmd)}{\Psi(\lmd)}\r)} & \text{if } \vep=1, \\[1pt]\\ 
\sqrt{\dfrac{\Phi(\lmd)}{\Psi(\lmd)} 
	+ \dfrac{2}{|1-\vep|\pi\lmd}\lf(\dfrac{\Phi(\lmd)}{\Psi(\lmd)}-1\r)} & \text{if } \vep\neq1, 
 \end{cases}
\end{align*}
then we find that 
\begin{align}\label{F4}
\begin{split}
\lf.\sum_{i=1}^{n-1}\frac{d^2L_i}{ds^2}\r|_{s=0} 
	&\leq -\frac{(n-1)HL\Psi(\lmd)}{2} \lf( 1- \frac{\pi^2}{HL^2} \widetilde{D}(\lmd)^2\r)\\
	&\quad + \Phi(\lmd)\int_0^L \big((n-1)H - \Ric_{N-} \big)_+ \,dt.
\end{split}
\end{align}
Substituting \eqref{L-com} for \eqref{F4} gives
\begin{align*}
\lf.\sum_{i=1}^{n-1}\frac{d^2L_i}{ds^2}\r|_{s=0} 
&\leq -\frac{(n-1)HL\Psi(\lmd) }{2}\lf(1 - \frac{\pi^2}{HL^2}\widetilde{D}(\lmd)^2 \r)\\
&\quad + 2r\,\lf(\frac{2b^2D}{a^2r}\r)^{\frac{1}{c}+1}\frac{\Phi(\lmd)}{\mu_f (M)} \int_M \big((n-1)H - \Ric_{N-}\big)_+\,d\mu_f .
\end{align*}
Choose $T = T(\eta)>2$ such that 
\begin{equation}\label{eta-inquality}
\frac{1}{1-\frac{2}{T}} \leq \frac{\pi + \eta}{\pi +\frac{\eta}{2}}
\end{equation}
and let $r = D/T$. By the assumption $y_i \in \overline{B}(p_i,r)$ for $i=1,2$ together with the triangle inequality, we find that
\begin{equation}\label{tri-inquality}
L=d(y_1,y_2) \geq d(p_1,p_2) - 2r =D\lf(1-\frac{2}{T}\r).
\end{equation}
Therefore we have
\begin{align*}
\lf.\sum_{i=1}^{n-1}\frac{d^2L_i}{ds^2}\r|_{s=0} 
&\leq -\frac{(n-1)HL\Psi(\lmd)}{2}\lf(1 - \frac{\pi^2}{HL^2}\widetilde{D}(\lmd)^2 \r)\\
&\quad +2\,\lf(\frac{2b^2}{a^2}\r)^{\frac{1}{c}+1}
T^{\frac{1}{c}}\frac{L}{1-\frac{2}{T}}\frac{\Phi(\lmd)}{\mu_f (M)} \int_M \big((n-1)H - \Ric_{N-}\big)_+\,d\mu_f .
\end{align*}
We consider the limit as
\begin{align}\label{lmd-0}
\lmd\to\lmd_0=
\begin{cases}
0 & \text{if } N\in [n,\infty],\\ 
\dfrac{1}{1-\vep}\lf(1-\sqrt{\dfrac{N-1}{N-n}}\r)& \text{if } N\in(-\infty, 1].
 \end{cases}
\end{align}
Then we have $\widetilde{D}(\lmd) \to \widetilde{D}(n,N,\vep,a,b)$ as $\lmd\to\lmd_0$. We set
\begin{align*}
\delta_1(n,N,\vep,a,b,H,\eta)
&= \frac{H(n-1)(1-\frac{2}{T})}{2^{\frac{1}{c}+3}T^{\frac{1}{c}}}\lf(1-\frac{\pi^2}{\lf(\pi +\frac{\eta}{2}\r)^2}\r) \lf(\frac{a}{b}\r)^{\frac{2}{c}+2+\lmd_0}.
\end{align*}
If we assume that
\begin{equation*}
\frac{1}{\mu_f (M)} \int_M \big((n-1)H - \Ric_{N-}\big)_+ \,d\mu_f \leq \delta_1(n,N,\vep,a,b,H,\eta),
\end{equation*} 
then we get
\begin{align*}
 \lf.\sum_{i=1}^{n-1}\frac{d^2L_i}{ds^2}\r|_{s=0}
&\leq  -\frac{(n-1)HLa^{\lmd_0} }{2}\lf(1 - \frac{\pi^2}{HL^2}\widetilde{D}^2 \r) 
	+ \frac{(n-1)HLa^{\lmd_0} }{2}\lf(1-\frac{\pi^2}{\lf(\pi+\frac{\eta}{2}\r)^2}\r) \\
&= -\frac{(n-1)H\pi^2a^{\lmd_0} }{2L\lf(\pi+\frac{\eta}{2}\r)^2}\lf(L^2 -\frac{1}{H}\lf(\pi+\frac{\eta}{2}\r)^2\widetilde{D}^2\r).
\end{align*}
Since $\gamma$ is a minimal geodesic, it follows
\begin{equation*}
 \lf.\sum_{i=1}^{n-1}\frac{d^2L_i}{ds^2}\r|_{s=0}\geq0.
\end{equation*}
Therefore we obtain
\begin{equation*}
L \leq \frac{1}{\sqrt{H}}\lf(\pi+\frac{\eta}{2}\r) \widetilde{D}(n,N,\vep,a,b). 
\end{equation*}
We observe from \eqref{eta-inquality}, \eqref{tri-inquality} that
\begin{equation*}
D
\leq \frac{L}{1-\frac{2}{T}} 
\leq \frac{\pi+\eta}{\sqrt{H}} \widetilde{D}(n,N,\vep,a,b).
\end{equation*}
This completes the proof of the theorem.
\end{apf}
\begin{rem}
In Theorem~\ref{Cpt-Main}, when $\vep\neq1$, if we take $a',b'$ such that
\begin{equation*}
\frac{b'}{a'}
=\exp\lf(\frac{2|1-\vep|}{n-1} (\sup f-\inf f)\r),
\end{equation*}
then $\delta_1(n,N,\vep,a,b,H,\eta) \leq \delta_1(n,N,\vep,a',b',H,\eta)$ holds. This implies that the assumption of Theorem~\ref{Cpt-Main} also holds when we replace $a,b$ by $a',b'$. We find that
\begin{align*}
\widetilde{D}(n,N,\vep,a',b')=
\begin{cases}
\sqrt{1+\dfrac{4(\sup f-\inf f) }{(n-1)\pi}} &\quad \text{if } N\in[n, \infty], \\[1pt]\\ 
\sqrt{\lmd_1 + \dfrac{2(n-N)(\lmd_1-1)}{(n-1)\pi} \left(1+\sqrt{\dfrac{N-1}{N-n}}\right)}&\quad \text{if } N\in(-\infty, 1],
 \end{cases}
\end{align*}
where
\begin{equation*}
\lmd_1=\exp\left(\dfrac{2(\sup f-\inf f)}{n-1}\left(1-\sqrt{\frac{N-1}{N-n}}\right)\right).
\end{equation*}
\end{rem}
%
\section{Proof of Theorem~\ref{Cmp}}\label{M-complete}
%
%
%
%
We provide the following lemma to prove Theorem~\ref{Cmp}.
\begin{lem}\label{Lem-cmp}
Let $n,N,\varepsilon, K, a, b,H,R$ and $\eta$ as in {\rm{Theorem~\ref{Cmp}}}. Assume $R, \eta$ satisfy
\begin{equation}\label{R-eta-cond}
R>\frac{\pi}{\sqrt{H}}\widetilde{D}(n,N,\vep,a,b) \quad\text{and}\quad
0<\eta < \eta_*(H,R,\widetilde{D})= \frac{4}{7}\lf(\frac{R\sqrt{H}}{\widetilde{D}}-\pi\r).
\end{equation}
There exists a positive constant $\delta_2(n,N,K,\vep,a,b,H,R,\eta)$ such that if an $n$-dimensional complete weighted manifold $(M,g,\mu_f )$  satisfies the $(N,K,\vep,a,b)$-condition and
\begin{equation}\label{Sec-ICB}
\sup_{p\in M} \frac{1}{\mu_f \big(B(p,R)\big)} \int_{B(p,R)} \big((n-1)H - \Ric_{N-}\big)_+ \,d\mu_f \leq \delta_2(n,N,K,\vep,a,b,H,R,\eta)
\end{equation} 
then $M$ is compact and
\begin{equation*}
\diam(M) \leq \frac{\pi+\eta}{\sqrt{H}}\widetilde{D}(n,N,\vep,a,b).
\end{equation*}
\end{lem}
\begin{pf}
The proof goes by contradiction, that is, there exist points $p_1,q\in M$ such that the distance from $p_1$ to $q$ is greater than $(\pi+\eta)\widetilde{D}/\sqrt{H}$. Then there exists $p_2 \in M$ such that $p_2$ lies in a unit minimal geodesic from $p_1$ to $q$ and
\begin{equation*}
\frac{\pi +\eta}{\sqrt{H}}\widetilde{D}
< d(p_1,p_2)
< R-\frac{3\eta\widetilde{D}}{4\sqrt{H}}.
\end{equation*}

First we set $W = B\lf(p_1, R\r)$ for $p_1 \in M$ and $r=\eta\widetilde{D}/4\sqrt{H}$ for $\eta<\eta_*(H,R,\widetilde{D})$. We put $A_i=B(p_i,r) \subset W$ for $i=1,2$. The triangle inequality that for $y_i\in A_i$ with $i=1,2$ yields
\begin{equation*}
d(y_1,y_2) \leq d(y_1,p_1)+d(p_1,p_2)+d(p_2,y_2)<R-r. 
\end{equation*} 
On the other hand, the distance from $y_1 \in A_1$ to the boundary of $W$ is greater than $R-r$. This means that all unit minimal geodesics from $y_1\in A_1$ to $y_2 \in A_2$ lie in $W$. Using Theorem~\ref{Seg-Ran}, we see that
\begin{align*}
&\quad \inf_{(y_1,y_2) \in \overline{A}_1 \times \overline{A}_2}
\int_0^{d(y_1,y_2)} \big((n-1)H - \Ric_{N-}\big)_+(\gamma_{y_1,y_2})\,ds \\
&\leq 2rC(n,c,K)\lf(\frac{1}{\mu_f (A_1)} + \frac{1}{\mu_f (A_2)}\r)\int_{B(p,R)} \big((n-1)H - \Ric_{N-}\big)_+\,d\mu_f .
\end{align*}
We set
\begin{equation*}
v_{cK,a}(R)=\int_0^{R/a} \sn_{cK}(\tau)^{1/c}\,d\tau.
\end{equation*} 
Theorem~\ref{BCT} implies
\begin{equation*}
\dfrac{\mu_f \big(B(p_1,R)\big)}{\mu_f \big(B(p_1,r)\big)}
\leq \frac{b}{a}\cdot\frac{v_{cK,a}(R)}{v_{cK,b}(r)}
\end{equation*}
and Lemma~\ref{Inc-sn} gives
\begin{equation*}
C(n,c,K)
\leq \widetilde{C} =\frac{\sn_{cK} \big((R-r)/a\big)^{1/c}} {\sn_{cK} \big(( R-r)/b\big)^{1/c}}.
\end{equation*}
We observe that
\begin{align}\label{L-inf}
\begin{split}
&\quad \inf_{(y_1,y_2) \in \overline{A}_1 \times \overline{A}_2}
\int_0^{d(y_1,y_2)} \big((n-1)H - \Ric_{N-}\big)_+(\gamma_{y_1,y_2})\,ds \\
&\leq 2rC(n,c,K)\frac{b}{a}\lf(
\frac{v_{cK,a}\lf(R\r)}{v_{cK,b}(r)}\frac{1}{\mu_f \big(B(p_1,R)\big)} +\frac{v_{cK,a}\lf(2R\r)}{v_{cK,b}(r)} \frac{1}{\mu_f \big(B(p_2, 2R)\big)}\r)\\
&\quad\times \int_{B(p_1,R)} \big((n-1)H - \Ric_{N-}\big)_+ \,d\mu_f \\
&\leq 2r\widetilde{C}\frac{b}{a}\lf(
\frac{v_{cK,a}\lf(R\r) + v_{cK,a}\lf(2R\r)}{v_{cK,b}(r)}
\r) \frac{1}{\mu_f \big(B(p_1,R)\big)} \int_{B(p_1,R)} \big((n-1)H - \Ric_{N-}\big)_+ \,d\mu_f ,
\end{split}
\end{align}
where we use $B(p_1, R)\subset B(p_2, 2R)$ in the second inequality.

We find a unit speed minimal geodesic $\gamma$ from $y_1 \in \overline{A}_1$ to $y_2 \in \overline{A}_2$ that attains the infimum of~\eqref{L-inf}. 
We put $L=d(y_1,y_2)$. With the same argument of Theorem~\ref{Cpt-Main}, we utilize \eqref{F4} and take the limit as $\lmd \to \lmd_0$ again, we see that
\begin{align*}
\lf.\sum_{i=1}^{n-1}\frac{d^2L_i}{ds^2}\r|_{s=0} 
&\leq  -\frac{(n-1)HLa^{\lmd_0} }{2}\lf(1 - \frac{\pi^2}{HL^2}\widetilde{D}^2 \r) 
 + \frac{\eta\widetilde{D}}{2\sqrt{H}}\frac{b^{\lmd_0+1}\widetilde{C}}{a}
\lf(
\frac{v_{cK,a}\lf(R\r) + v_{cK,a}\lf(2R\r)}{v_{cK,b}\lf(\frac{\eta\widetilde{D}}{4\sqrt{H}}\r)}
\r)\\
&\quad\times \frac{1}{\mu_f \big(B(p,R)\big)}\int_{B(p,R)} \big((n-1)H - \Ric_{N-}\big)_+\,d\mu_f .
\end{align*}
We set
\begin{align}\label{Sec-ICC}
\begin{split}
 &\quad \delta_2(n,N,K,\vep,a,b,H,R,\eta)\\
&= \frac{H(n-1)(\pi +\frac{\eta}{2})}{\eta \widetilde{C}}
\frac{v_{cK,b}\lf(\frac{\eta\widetilde{D}}{4\sqrt{H}}\r)}{v_{cK,a}\lf(R\r) + v_{cK,a}\lf(2R\r)}
\lf(
1-\frac{\pi^2}{\lf(\pi +\frac{\eta}{2}\r)^2}
\r)\lf(\frac{a}{b}\r)^{\lmd_0+1}.
\end{split}
\end{align}
If we assume that
\begin{equation*}
\sup_{p\in M} \frac{1}{\mu_f \big(B(p,R)\big)} \int_{B(p,R)} \big((n-1)H - \Ric_{N-}\big)_+ \,d\mu_f \leq \delta_2(n,N,K,\vep,a,b,H,R,\eta),
\end{equation*} 
then we have
\begin{align*}
\lf.\sum_{i=1}^{n-1}\frac{d^2L_i}{ds^2}\r|_{s=0}
&\leq -\frac{(n-1)HLa^{\lmd_0}}{2}\lf(1-\frac{\pi^2}{HL^2}\widetilde{D}^2 \r) \\
&\quad+ \frac{(n-1)(\pi +\frac{\eta}{2})\sqrt{H}a^{\lmd_0}\widetilde{D}}{2}
\lf(
1-\frac{\pi^2}{(\pi +\frac{\eta}{2})^2}
\r)\\
&= -\frac{n-1}{2}
\lf(L - \frac{1}{\sqrt{H}}\lf(\pi +\frac{\eta}{2}\r)\widetilde{D}\r)\lf(
Ha^{\lmd_0} + \frac{(n-1)\pi^2\widetilde{D}a^{\lmd_0}\sqrt{H}}{2L(\pi +\frac{\eta}{2})}
\r).
\end{align*}
Since $\gamma$ is a minimal geodesic, it follows
\begin{equation*}
 \lf.\sum_{i=1}^{n-1}\frac{d^2L_i}{ds^2}\r|_{s=0}\geq0.
\end{equation*}
Therefore we obtain
\begin{equation*}
L \leq \frac{1}{\sqrt{H}}\lf(\pi +\frac{\eta}{2}\r)\widetilde{D}.
\end{equation*}
By the triangle inequality, we have 
\begin{equation}\label{contra}
d(p_1,p_2) \leq L +2r \leq \frac{\pi +\eta}{\sqrt{H}}\widetilde{D}.
\end{equation}
On the other hand, we assumed that
\begin{equation*}
\frac{\pi +\eta}{\sqrt{H}}\widetilde{D}
< d(p_1,p_2)
< R-\frac{3\eta\widetilde{D}}{4\sqrt{H}},
\end{equation*}
but by \eqref{contra}, no geodesic starting from $p_1 \in M$ of a length greater than $(\pi+\eta)\widetilde{D}/\sqrt{H}$ can be minimal. 
This completes the proof of the lemma.
\end{pf}
Finally, we show Theorem~\ref{Cmp}.
\begin{apf}{Theorem~\ref{Cmp}}
We divide the cases into three parts:
\begin{enumerate}
\item \label{R-ls} $R \leq \pi\widetilde{D}(n,N,\vep,a,b)/\sqrt{H}$,
\item \label{R-bg} $R>\pi\widetilde{D}(n,N,\vep,a,b)/\sqrt{H}$ and $\eta \geq \eta_*(H,R,\widetilde{D})$,
\item \label{R-lmm}$R>\pi\widetilde{D}(n,N,\vep,a,b)/\sqrt{H}$ and $\eta < \eta_*(H,R,\widetilde{D})$.
\end{enumerate}
The case \eqref{R-lmm} is already discussed in Lemma~\ref{Lem-cmp}. Thus it is enough to consider the cases \eqref{R-ls} and \eqref{R-bg}. 

Let $R'=R'(\eta) > \pi\widetilde{D}/\sqrt{H}$ be fixed such that $\eta<\eta_*(H,R',\widetilde{D})$. For each $p\in M$, we consider the discrete subset $\{x_i\} \subset B(p,R')$ for $i=1,\ldots ,T_2$ such that 
\begin{equation*}
B(p,R') \subset \bigcup_{i=1}^{T_2} B(x_i, R)
\end{equation*}
 and $d(x_i,x_j) >R$ for $i\neq j$, where $T_2$ is the maximal number of the $R$-discrete net of $B(p,R')$ (see~\cite{SHT}*{Definition~3.1} for the definition of the $R$-discrete net). We now claim:
\begin{claim}Let $(M,g,\mu_f )$ be an $n$-dimensional complete weighted manifold satisfying the $(N,K,\vep,a,b)$-condition. Set $T_2>0$ as above. Then it follows 
\begin{equation*}
T_2 \leq \frac{b}{a}\cdot\frac{\vb (2R'+R)}{\va\lf(R/2\r)}.
\end{equation*}
\end{claim}
\begin{cpf}
Take $i_0 \in \{\,1,\ldots,T_2\,\}$ such that
\begin{equation*}
\mu_f\big(B(x_{i_0},R/2)\big)=\min_{1\leq i\leq T_2} \mu_f\big(B(x_i,R/2)\big).
\end{equation*}
 Since $B(x_i,R/2) \cap B(x_j,R/2) = \emptyset$ for any $i\neq j$ and $B(p,R'+R) \subset B(x_{i_0},2R'+R)$ hold, we have  
\begin{align*}
1
&= \frac{\mu_f\big(B(p,R'+R)\big)}{\mu_f\big(B(p,R'+R)\big)} \\
&\leq \frac{\mu_f\big(B(p,R'+R)\big)}{\mu_f\big(B(p,R'+R) \cap \bigcup_{i=1}^{T_2}B(x_i,R/2)\big)} \\
&\leq \frac{\mu_f\big(B(p,R'+R)\big)}{T_2\cdot\mu_f\big(B(x_{i_0},R/2)\big)} \\
&\leq \frac{\mu_f\big(B(x_{i_0},2R'+R)\big)}{T_2\cdot\mu_f\big(B(x_{i_0},R/2)\big)}.
\end{align*}
By using Theorem~\ref{BCT}, we obtain
\begin{equation*}
T_2\leq \frac{b}{a}\cdot\frac{\vb (2R'+R)}{\va(R/2)}.
\end{equation*}
\end{cpf}
Then we find that, for any $z\in M$,  
\begin{align*}
&\quad \frac{1}{\mu_f \big(B(z,R')\big)} \int_{B(z,R')} \big((n-1)H - \Ric_{N-}\big)_+\,d\mu_f \\
&\leq \frac{T_2}{\mu_f \big(B(z,R')\big)}\sup_{x_i \in B(z,R')} \int_{B(x_i,R)} \big((n-1)H - \Ric_{N-}\big)_+\,d\mu_f \\
&\leq T_2 \frac{b}{a} \frac{1}{\mu_f \big(B(z,R' + R)\big)}
\frac{v_{cK,b}(R'+R)}{v_{cK,a}(R')}\\
&\quad \times \sup_{x_i\in B(z,R')} \int_{B(x_i,R)} \big((n-1)H - \Ric_{N-}\big)_+\,d\mu_f \\
&\leq \frac{b^2}{a^2}\frac{\vb (2R'+R)}{\va(R/2)}
\frac{v_{cK,b}(R'+R)}{v_{cK,a}(R')}\\
&\quad\times\sup_{x_i \in B(z,R')} \frac{1}{\mu_f \big(B(x_i,R)\big)}\int_{B(x_i,R)} \big((n-1)H - \Ric_{N-}\big)_+\,d\mu_f ,
\end{align*}
where we use the estimate of $T_2$ and $B(x_i,R) \subset B(z,R' + R)$ for $1\leq i\leq T_2$ in the third inequality. This provides that
\begin{align*}
&\quad \sup_{x\in M} \frac{1}{\mu_f \big(B(x,R')\big)} \int_{B(x,R')} \big((n-1)H - \Ric_{N-}\big)_+\,d\mu_f \\
&\leq \frac{b^2}{a^2}\frac{\vb (2R'+R)}{\va\lf(R/2\r)}\frac{v_{cK,b}(R'+R)}{v_{cK,a}(R')}
	\sup_{x\in M} \frac{1}{\mu_f\big(B(x,R)\big)}\int_{B(x,R)} \big((n-1)H - \Ric_{N-}\big)_+\,d\mu_f.
\end{align*}
We set
\begin{align*}
&\quad \delta'_2(n,N,K,\vep,a,b,H,R,\eta)\\
&= \frac{a^2}{b^2}\frac{\va(R/2)}{\vb (2R'(\eta)+R)} \frac{v_{cK,a}\big(R'(\eta)\big)}{v_{cK,b}\big(R'(\eta)+R\big)}\delta_2(n,N,K,\vep,a,b,H,R'(\eta),\eta).
\end{align*}
If we assume that
\begin{equation*}
\sup_{x\in M} \frac{1}{\mu_f\big(B(x,R)\big)}\int_{B(x,R)}\big((n-1)H - \Ric_{N-}\big)_+ \,d\mu_f \leq \delta'_2(n,N,K,\vep,a,b,H,R,\eta),
\end{equation*} 
we obtain
\begin{equation*}
\sup_{x\in M} \frac{1}{\mu_f \big(B(x,R')\big)} \int_{B(x,R')} \big((n-1)H - \Ric_{N-}\big)_+\,d\mu_f 
\leq \delta_2(n,N,K,\vep,a,b,H,R'(\eta),\eta).
\end{equation*}
By Lemma~\ref{Lem-cmp}, this completes the proof of the cases \eqref{R-ls}, \eqref{R-bg} and the theorem.
\end{apf}
We prepare the following lemma to provide a slightly weak result compared with Lemma~\ref{Lem-cmp}. 
\begin{lem}\label{Dlt-Inc}
Let $\delta_2=\delta_2(n,N,K,\vep,a,b,H,R,\eta)$ be given as in \eqref{Sec-ICC}. Then $\delta_2$ is strictly increasing in $\eta\in (0,\eta_*)$. 
\end{lem}

\begin{pf}
Set
\begin{equation}\label{geta}
G(\eta)= \lf(\frac{1}{2} +\frac{\pi}{\eta}\r)\int_0^{\eta\widetilde{D}/4b\sqrt{H}} \sn_{cK}(t)^{\frac{1}{c}}\, dt.
\end{equation}
Note that the monotonicity of $\sn_{cK}$. We observe that
\begin{align*}
G'(\eta)
&= \frac{\widetilde{D}}{4b\sqrt{H}}\lf(\frac{1}{2} +\frac{\pi}{\eta}\r) \lf(\sn_{cK}\lf(\frac{\eta\widetilde{D}}{4b\sqrt{H}}\r)\r)^{\frac{1}{c}}
		-\frac{\pi}{\eta^2}\int_0^{\eta\widetilde{D}/4b\sqrt{H}} \sn_{cK}(t)^{\frac{1}{c}}\, dt \\
&> \frac{\widetilde{D}}{4b\sqrt{H}}\lf(\frac{1}{2} +\frac{\pi}{\eta}\r) \lf(\sn_{cK}\lf(\frac{\eta\widetilde{D}}{4b\sqrt{H}}\r)\r)^{\frac{1}{c}}
		-\frac{\pi}{\eta^2} \frac{\eta\widetilde{D}}{4b\sqrt{H}}\lf(\sn_{cK}\lf(\frac{\eta\widetilde{D}}{4b\sqrt{H}}\r)\r)^{\frac{1}{c}}\\
&= \frac{\widetilde{D}}{8b\sqrt{H}} \lf(\sn_{cK}\lf(\frac{\eta\widetilde{D}}{4b\sqrt{H}}\r)\r)^{\frac{1}{c}}\\
&>0.
\end{align*}
Thus $G(\eta)$ is strictly increasing on $(0,\eta_*)$. Moreover
\begin{equation}\label{etae}
1-\frac{\pi^2}{\lf(\pi+\frac{\eta}{2}\r)^2}
\end{equation}
and, by Claim~\ref{sn-clm}, 
\begin{equation}\label{etaee}
\frac{\sn_{cK} \Big(\frac{1}{b}\lf( R-\frac{\eta\widetilde{D}}{4\sqrt{H}}\r)\Big)^{\frac{1}{c}}}
{\sn_{cK} \Big(\frac{1}{a}\lf( R-\frac{\eta\widetilde{D}}{4\sqrt{H}}\r)\Big)^{\frac{1}{c}}}
\end{equation}
are also strictly increasing with respect to $\eta$. Since $\delta_2$ was given by multiplication of \eqref{geta},~\eqref{etae},~\eqref{etaee} and a positive constant independent of $\eta$, it follows that $\delta_2$ is strictly increasing in $\eta\in (0,\eta_*)$.
\end{pf}
We take the limit of \eqref{Sec-ICC} as $\eta \to \eta_*$, where $\eta_*$ is given in \eqref{R-eta-cond}. Then we have a diameter estimate of a weighted manifold. 
\begin{cor}\label{R-EST}
Let $n,N,\varepsilon, K, a, b, H,R$ and $\eta_*$ as in {\rm{Lemma~\ref{Lem-cmp}}}. 
There exists a positive constant $ \delta_2(n,N,K,\vep,a,b,H,R,\eta_*)$ such that if an $n$-dimensional complete weighted manifold $(M,g,\mu_f )$ satisfies the $(N,K,\vep,a,b)$-condition and
\begin{equation}\label{R-ICB}
\sup_{p\in M} \frac{1}{\mu_f \big(B(p,R)\big)} \int_{B(p,R)} \big((n-1)H - \Ric_{N-}\big)_+ \,d\mu_f < \delta_2(n,N,K,\vep,a,b,H,R,\eta_*)
\end{equation} 
then $M$ is compact and $\diam(M) < R$.
\end{cor}
\begin{pf}
By Lemma~\ref{Dlt-Inc}, if \eqref{R-ICB} holds, then there exists $\eta>0$ such that $\eta <\eta_*(H,R,\widetilde{D})$ satisfying \eqref{Sec-ICB}. Hence Lemma~\ref{Lem-cmp} implies
\begin{equation*}
\diam(M) \leq \frac{\pi+\eta}{\sqrt{H}}\widetilde{D}(n,N,\vep,a,b) <R.
\end{equation*}
\end{pf}
%
\section{Fundamental group under integral curvature bound and $\vep$-range}\label{FGB}
%
 We shall show a finiteness of  the fundamental group of $M$. Compared with Corollary~\ref{R-EST}, we need to assume a slightly strong condition about integral curvature bound.
\begin{cor}\label{FGB-lemma}
Let $(M,g,\mu_f )$ be an $n$-dimensional complete weighted manifold satisfying the $(N,K,\vep,a,b)$-condition. If there exists $H>0$, $R>\pi\tilde{D}(n,N,\vep,a,b)/\sqrt{H}$ such that
\begin{equation*}
\sup_{p\in M} \frac{1}{\mu_f \big(B(p,R)\big)} \int_{B(p,R)} \big((n-1)H - \Ric_{N-}\big)_+ \,d\mu_f 
< \frac{a}{b}\frac{v_{cK,b}(R)}{v_{cK,a}(3R)}\delta_2,
\end{equation*} 
where $\delta_2=\delta_2(n,N,K,\vep,a,b,H,R,\eta_*)$ is given in \eqref{Sec-ICC}, then the universal cover of $M$ is compact, and hence $\pi_1(M)$ is a finite group.
\end{cor}

\begin{pf}
We find that, by $3R/a \geq R/b$, 
\begin{equation*}
\frac{a}{b}\frac{v_{cK,b}(R)}{v_{cK,a}(3R)}<1.
\end{equation*}
Since $(M,g,\mu_f )$ satisfies the assumption of Corollary~\ref{R-EST}, we have $\diam(M) <R$.

We denote by $k : \widetilde{M} \to M$ the universal Riemannian covering. Set $\tilde{g} = k^*g$, $\tilde{f}=f\circ k$. Fix $\tilde{x}\in \widetilde{M}$. 
Let $F$ be the fundamental domain of $k$ that contains $\tilde{x}\in \widetilde{M}$. We find that
\begin{equation*}
\mu_f(M)= \mu_{\tilde{f}}(F),
\qquad
\widetilde{M} =\bigcup_{\alpha \in \Gamma} \alpha F,
\end{equation*}
where $\Gamma$ is the deck transformation group of $\widetilde{M}$.
We set
\begin{equation*}
T= \inf \lf\{\# \Gamma_0 \,\middle|\, B(\tilde{x}, R) \subset \bigcup_{\alpha \in \Gamma_0} \alpha F\r\}.
\end{equation*}
Fix $w \in B(\tilde{x}, R) \cap \alpha F$ and $z\in \alpha F$ for each $\alpha\in\Gamma_0$. Since $d(y,z) < 2R$ for all $y \in \alpha F$, we see that 
\begin{equation*}
d(\tilde{x}, z) \leq d(\tilde{x},w)+d(w,z)<3R.
\end{equation*}
Hence this implies
\begin{equation*}
\bigcup_{i=1}^T \alpha_i F \subset B(\tilde{x},3R)
\end{equation*}
and we find that 
\begin{equation}\label{TVV}
T\cdot \mu_f (M)= T\cdot \mu_{\tilde{f}}(F)\leq \mu_{\tilde{f}} \big(B(\tilde{x},3R)\big).
\end{equation}
Then it follows from \eqref{TVV} that
\begin{align*}
&\quad \frac{1}{\mu_{\tilde{f}}  \big(B(\tilde{x},R)\big)} \int_{B(\tilde{x},R)} \big((n-1)H - \Ric_{N-}\big)_+\,d\mu_{\tilde{f}} \\
&\leq \frac{T}{\mu_{\tilde{f}} \big(B(\tilde{x},R)\big)} \int_{F} \big((n-1)H - \Ric_{N-}\big)_+\,d\mu_{\tilde{f}} \\
&= \frac{T}{\mu_{\tilde{f}} \big(B(\tilde{x},R)\big)} \int_{M} \big((n-1)H - \Ric_{N-}\big)_+\,d\mu_f \\
&\leq \frac{b}{a}
\frac{v_{cK,b}(3R)}{v_{cK,a}(R)}\frac{T}{\mu_{\tilde{f}} \big(B(\tilde{x},3R)\big)}\int_{M} \big((n-1)H - \Ric_{N-}\big)_+\,d\mu_f \\
&\leq \frac{b}{a}
\frac{v_{cK,b}(3R)}{v_{cK,a}(R)}\frac{1}{\mu_f (M)}\int_{M} \big((n-1)H - \Ric_{N-}\big)_+\,d\mu_f .
\end{align*}
Thus we observe that
\begin{align*}
&\quad \sup_{\tilde{x}\in \widetilde{M}} \frac{1}{\mu_{\tilde{f}}  \big(B(\tilde{x},R)\big)} \int_{B(\tilde{x},R)} \big((n-1)H - \Ric_{N-}\big)_+\,d\mu_{\tilde{f}} \\
&\leq \frac{b}{a}
\frac{v_{cK,b}(3R)}{v_{cK,a}(R)} \frac{1}{\mu_f (M)}\int_M \big((n-1)H - \Ric_{N-}\big)_+\,d\mu_f \\
&<\delta_2.
\end{align*}
Since $(\widetilde{M},\tilde{g},\mu_{\tilde{f}})$ satisfies the assumption of Corollary~\ref{R-EST}, $\widetilde{M}$ is compact and $\pi_1(M)$ is a finite group. 
\end{pf}
Finally, we extend Corollary~\ref{FGB-lemma} to the cases that $R \leq \pi\widetilde{D}/\sqrt{H}$.
\begin{cor}
Let $n,N,\varepsilon, K, a, b, H$ and $R$ as in {\rm{Theorem~\ref{Cmp}}}. Set $R'>0$ satisfies $R'>\pi\widetilde{D}(n,N,\vep,a,b)/\sqrt{H}$. There exists a positive constant $\tilde{\delta}(n,N,K,\vep,a,b,H,R,R')$ such that if an $n$-dimensional complete weighted manifold $(M,g,\mu_f )$ satisfies the $(\vep,a,b,N,K)$-condition and
\begin{equation*}
\sup_{p\in M} \frac{1}{\mu_f \big(B(p,R)\big)} \int_{B(p,R)} \big((n-1)H - \Ric_{N-}\big)_+ \,d\mu_f 
< \tilde{\delta}(n,N,K,\vep,a,b,H,R,R')
\end{equation*} 
then the universal cover of $M$ is compact, and hence $\pi_1(M)$ is a finite group.
\end{cor}
\begin{pf}
Let $R \leq \pi\widetilde{D}/\sqrt{H}$.
With the same argument in the proof of Theorem~\ref{Cmp}, we find that
\begin{align*}
&\quad \sup_{x\in M} \frac{1}{\mu_f \big(B(x,R')\big)} \int_{B(x,R')} \big((n-1)H - \Ric_{N-}\big)_+\,d\mu_f \\
&\leq \frac{b^2}{a^2}\frac{\vb (2R'+R)}{\va(R/2)}\frac{v_{cK,b}(R'+R)}{v_{cK,a}(R')}
	\sup_{x\in M} \frac{1}{\mu_f \big(B(x,R)\big)} \int_{B(x,R)} \big((n-1)H - \Ric_{N-}\big)_+\,d\mu_f .
\end{align*}
If we assume that 
\begin{align*}
&\quad \sup_{x\in M} \frac{1}{\mu_f \big(B(x,R)\big)} \int_{B(x,R)} \big((n-1)H - \Ric_{N-}\big)_+\,d\mu_f\\
&<\frac{a^3}{b^3}\frac{\va(R/2)}{\vb (2R'+R)}\frac{v_{cK,a}(R')}{v_{cK,b}(R'+R)}
 \frac{v_{cK,b}(R')}{v_{cK,a}(3R')}\delta_2(n,N,K,\vep,a,b,H,R',\eta_*)
\end{align*}
holds, we obtain
\begin{align*}
&\quad \sup_{x\in M} \frac{1}{\mu_f  \big(B(x,R')\big)} \int_{B(x,R')} \big((n-1)H - \Ric_{N-}\big)_+\,d\mu_f \\
&<\frac{a}{b} \frac{v_{cK,b}(R')}{v_{cK,a}(3R')}\delta_2(n,N,K,\vep,a,b,H,R',\eta_*).
\end{align*}
Therefore, by Corollary~\ref{FGB-lemma}, this completes the proof of the corollary.
\end{pf}


\begin{bibdiv}
 \begin{biblist}
 %
%
%
%
%
\bib{CC}{article}{
   author={Cheeger, Jeff},
   author={Colding, Tobias H.},
   title={Lower bounds on Ricci curvature and the almost rigidity of warped
   products},
   journal={Ann. of Math. (2)},
   volume={144},
   date={1996},
   number={1},
   pages={189--237},
   issn={0003-486X},
}
\bib{MR3952638}{article}{
   author={Hwang, Seungsu},
   author={Lee, Sanghun},
   title={Integral curvature bounds and bounded diameter with Bakry-\'{E}mery
   Ricci tensor},
   journal={Differential Geom. Appl.},
   volume={66},
   date={2019},
   pages={42--51},
   issn={0926-2245},
} 
%
\bib{HL}{article}{
   author={Hwang, Seungsu},
   author={Lee, Sanghun},
   title={Erratum to: ``Integral curvature bounds and bounded diameter with
   Bakry-\'{E}mery Ricci tensor'' [Differ. Geom. Appl. 66 (2019) 42--51]},
   journal={Differential Geom. Appl.},
   volume={70},
   date={2020},
   pages={101627, 3},
   issn={0926-2245},
}
%
\bib{JM}{article}{
  author={Jaramillo, Maree},
   title={Fundamental groups of spaces with Bakry-Emery Ricci tensor bounded
  below},
   journal={J. Geom. Anal.},
   volume={25},
   date={2015},
   number={3},
   pages={1828--1858},
   issn={1050-6926},
}
%
\bib{KL}{arXiv}{
     title={New Laplacian comparison theorem and its applications to diffusion processes on Riemannian manifolds}, 
     author={Kuwae, Kazuhiro},
    author={Li, Xiang-Dong},
      year={2020},
      eprint={2001.00444},
      archivePrefix={arXiv}
}
\bib{KS}{article}{
   author={Kuwae, Kazuhiro},
   author={Sakurai, Yohei},
   title={Rigidity phenomena on lower $N$-weighted Ricci curvature bounds
   with $\varepsilon$-range for nonsymmetric Laplacian},
   journal={Illinois J. Math.},
   volume={65},
   date={2021},
   number={4},
   pages={847--868},
   issn={0019-2082},
}
\bib{LMO}{article}{
   title={Geometry of weighted Lorentz–Finsler manifolds I: singularity theorems},
  volume={104},
   ISSN={1469-7750},
   journal={Journal of the London Mathematical Society},
  publisher={Wiley},
   author={Lu, Yufeng},%
   author={Minguzzi, Ettore},
   author={Ohta, Shin‐ichi},
   year={2021},
   month={Jan},
	pages={362–393}
}
\bib{LMO2}{arXiv}{
      title={Comparison theorems on weighted Finsler manifolds and spacetimes with $\epsilon$-range}, 
      author={Lu, Yufeng},
      author={Minguzzi, Ettore},
      author={Ohta, Shin‐ichi},
      year={2020},
      eprint={2007.00219},
      archivePrefix={arXiv}
}
 \bib{SBM}{article}{
   author={Myers, S. B.},
   title={Riemannian manifolds with positive mean curvature},
   journal={Duke Math. J.},
   volume={8},
   date={1941},
   pages={401--404},
   issn={0012-7094},
}
\bib{ST}{book}{
   author={Sakai, Takashi},
   title={Riemannian geometry},
   series={Translations of Mathematical Monographs},
   volume={149},
   note={Translated from the 1992 Japanese original by the author},
   publisher={American Mathematical Society, Providence, RI},
   date={1996},
   pages={xiv+358},
   isbn={0-8218-0284-4},
}
%
\bib{SHT}{book}{
   author={Shioya, Takashi},
   title={Metric measure geometry},
   series={IRMA Lectures in Mathematics and Theoretical Physics},
   volume={25},
   publisher={EMS Publishing House, Z\"{u}rich},
   date={2016},
   pages={xi+182},
   isbn={978-3-03719-158-3},
}
%
\bib{SC}{article}{
   author={Sprouse, Chadwick},
   title={Integral curvature bounds and bounded diameter},
   journal={Comm. Anal. Geom.},
   volume={8},
   date={2000},
   number={3},
   pages={531--543},
   issn={1019-8385},
}
%
%
\bib{WY}{arXiv}{
      title={On the geometry of Riemannian manifolds with density}, 
     author={Wylie, William},
      author={Yeroshkin,Dmytro},
      year={2016},%
      eprint={1602.08000},
      archivePrefix={arXiv}
}
%
%
 \end{biblist}
\end{bibdiv}
\end{document}